\documentclass[reqno]{amsart}
\usepackage{amsfonts,amsmath,amssymb,color}

\usepackage[babel=true,kerning=true]{microtype}

\usepackage{amsthm,leqno}

\usepackage{amsfonts,amsmath,amssymb,enumerate}
\usepackage{graphicx}
\usepackage{tikz}
\usetikzlibrary{decorations.markings}
\usetikzlibrary{plotmarks}
\usetikzlibrary{patterns}

\newtheorem{theo}{Theorem}

\newtheorem*{lem}{Lemma}

\begin{document}

\title[A rigidity result for the first conformal eigenvalue]{On a rigidity result for the first conformal eigenvalue of the Laplacian}
\author{Romain Petrides} 
\address{Romain Petrides, Universit\'e de Lyon, CNRS UMR 5208, Universit\'e Lyon 1, Institut Camille Jordan, 43 bd du 11 novembre 1918, F-69622 Villeurbanne cedex, France.}
\email{romain.petrides@univ-lyon1.fr}
\date{}

\begin{abstract} Given $(M,g)$ a smooth compact Riemannian manifold without boundary of dimension $n\geq 3$, we consider the first conformal eigenvalue which is by definition the supremum of the first eigenvalue of the Laplacian among all metrics conformal to $g$ of volume $1$. We prove that it is always greater than $n\omega_n^{\frac{2}{n}}$, the value it takes in the conformal class of the round sphere, except if $(M,g)$ is conformally diffeomorphic to the standard sphere. 
\end{abstract}

\maketitle

\medskip Let $(M,g)$ be a smooth compact Riemannian manifold without boundary of dimension $n\ge 3$ and let us define the first conformal eigenvalue of $(M,g)$ by 
$$ \Lambda_{1}(M,[g]) = \sup_{\tilde{g}\in[g]} \lambda_1(M,\tilde{g})Vol_{\tilde{g}}(M)^{\frac{2}{n}} $$
where $\lambda_1(M,g)$ is the first nonzero eigenvalue of the Laplacian $\Delta_g = -div_g\left(\nabla\right)$ and $[g]$ is the conformal class of $g$. In this paper, we aim at proving a rigidity result concerning this first conformal eigenvalue.

The maximisation on conformal classes is natural because the scale invariant quantity supremum is infinite among all metrics \cite{COL1994} (except in dimension $2$, \cite{YAN1980}), while El Soufi and Ilias \cite{ELS1986} proved that it is always bounded among conformal metrics. Generalizing a result by Li and Yau \cite{LIY1982} in dimension $2$, they gave an explicit upper bound thanks to the $m$-conformal volume $V_c(m,M,[g])$ of $(M,[g])$
\begin{equation} \label{boundconfvol} \Lambda_1(M,[g]) \leq n V_c(m,M,[g])^{\frac{2}{n}} \end{equation}
These conformal invariants on the standard sphere $(\mathbb{S}^n,[can])$ satisfy, \cite{ELS1986}
\begin{equation}\label{hers} \Lambda_{1}(\mathbb{S}^n,[can]) = n\omega_n^{\frac{2}{n}} = n V_c(\mathbb{S}^n,[can])^{\frac{2}{n}} \end{equation}
and this value is achieved if and only if the metric is round. Here, $\omega_n$ denotes the volume of the standard $n$-sphere. Colbois and El Soufi \cite{COL2003} also proved that, for any compact Riemannian manifold $(M,g)$ of dimension $n\geq 3$ 
$$ \Lambda_{1}(M,[g]) \geq \Lambda_{1}(\mathbb{S}^n,[can]) \hskip.1cm.$$ 
We prove here that the case of equality characterizes the standard sphere~:

\begin{theo} \label{rig} Let $(M,g)$ be a compact Riemannian manifold without boundary of dimension $n\geq 3$. Then
$$ \Lambda_{1}(M,[g]) > \Lambda_{1}(\mathbb{S}^n,[can])$$
if $(M,[g])$ is not conformally diffeomorphic to $(\mathbb{S}^n,[can])$.
\end{theo}

This theorem answers the question raised in \cite{COL2004} and \cite{KON2010}. Note that a similar result was proved by the author in dimension $2$ (see \cite{PETR2013}). Note also that thanks to (\ref{boundconfvol}) and (\ref{hers}), the theorem implies
$$ V_c(m,M,[g]) > \omega_n = V_c(\mathbb{S}^n,[can]) $$
if $(M,[g])$ is not conformally diffeomorphic to $(\mathbb{S}^n,[can])$. This gives a positive answer to question 2 in \cite{LIY1982}.

\medskip In the rest of this paper, we prove the theorem. Based on the idea of Ledoux \cite{LED2000} and Druet \cite{DRU2002}, we start from a sharp Sobolev inequality in dimensions $n \geq 3$ (see \cite{HEB1992,DRU2002,DRH2002}) which possesses extremal functions. These extremal functions give natural metrics $\tilde{g}\in[g]$ with $Vol_{\tilde{g}}(M) =1$ and $\lambda_1(\tilde{g}) \geq n\omega_n^{\frac{2}{n}}$.  As in dimension $2$, see \cite{PETR2013}, we deal with the degeneracy consequences of the hypothesis $\lambda_1(\tilde{g}) = n\omega_n^{\frac{2}{n}}$.

\medskip Let $(M,g)$ be a smooth compact Riemannian manifold of dimension $n\geq 3$ with $Vol_g(M) = 1$, which is not conformally diffeomorphic to the standard sphere. For an integer $m\geq 1$, let $h\in \mathcal{C}^{m}(M)$. We let  $J_{g,h}$ be the functional defined for $u\in W^{1,2}(M)\setminus \{0\}$ by 
\begin{equation} \label{defJ} J_{g,h}(u) = \frac{ \int_M \left|\nabla u\right|^2_g dv_g + \int_M h u^2 dv_g - K_n^{-2}\left(\int_M \left|u\right|^{2^{*}} dv_g \right)^{\frac{2}{2^{*}}} }{\int_M u^2 dv_g} \end{equation}
where 
\begin{equation} \label{defKn} 
K_n = \frac{2}{\sqrt{n(n-2)}}\omega_n^{-\frac{1}{n}}\end{equation}
is the sharp constant for the Sobolev inequality induced by the critical Sobolev embedding $W^{1,2}_0 \subset L^{2^{*}}$ for bounded domains of $\mathbb{R}^n$, with $2^{*} = \frac{2n}{n-2}$. Hebey and Vaugon proved in \cite{HEB1992} that 
\begin{equation} \label{defalpha} -\alpha(g,h) = \inf_{u\in W^{1,2}(M)\setminus \{0\}} J_{g,h}(u)
\end{equation}
is finite. Note that $J_{g,h}$ is scale invariant.

We will assume in the following that up to a conformal change, $g$ is a metric in $[g]$ with volume $1$ which has a constant scalar curvature $S_g$. Since $M$ is not conformally diffeomorphic to the standard sphere, by the resolution of the Yamabe problem by Aubin \cite{AUB1976} and Schoen \cite{SCH1984}, it satisfies
\begin{equation}\label{conditionmu} \mu(M,g) < K_n^{-2} \end{equation}
where $\mu(M,g)$ is the Yamabe invariant of $(M,[g])$. Let $V$ be an open neighbourhood of $\frac{n-2}{4(n-1)} S_g$ in $\mathcal{C}^m(M)$ such that
\begin{equation} \label{conddruet} \forall h\in V, \left\|h - \frac{n-2}{4(n-1)} S_g\right\|_{\infty} \leq \frac{1}{2} \left(K_n^{-2} - \mu(M,g) \right) \hskip.1cm. \end{equation}
Let $s\geq 0$ be such that $s+2> \frac{n}{2}$ and $m\geq s+2$. By the Sobolev embedding $W^{s+2,2}\hookrightarrow \mathcal{C}^0$, the subset $W^{s+2,2}_+$ of positive functions of $W^{s+2,2}$ is open. We define
$$\begin{array}{ccccc}
F & : & W^{s+2,2}_+ \times \mathbb{R} \times V & \longrightarrow & W^{s,2} \\
& & (u,\beta,h) & \longmapsto & \Delta_g u + (h + \beta) u - K_n^{-2} u^{2^{*}-1} \\
\end{array}$$
which is well defined because of the Sobolev algebra property of $W^{s+2,2}$ and $F$ is a $\mathcal{C}^{\infty}$ map. By a result of Druet \cite{DRU2002}, thanks to (\ref{conditionmu}) and (\ref{conddruet}), for any $h \in V$, the functional $J_{g,h}$ attains its infimum. Let $u \in W^{1,2}(M)$ be such that $J_{g,h}(u) = - \alpha(g,h)$. Up to replace $u$ by $\left|u\right|$ and up to normalize, we can take $u \geq 0$ and $\int_M u^{2^{*}}dv_g = 1$. Then, $u$ satisfies the Euler-Lagrange equation
\begin{equation}\label{eullageq} F(u,\alpha(g,h),h) = \Delta_g u + (h + \alpha(g,h)) u - K_n^{-2} u^{2^{*}-1} = 0 \end{equation}
where, by elliptic regularity theory, $u\in \mathcal{C}^{m+2}$ and, by the maximum principle, $u>0$. 

Let $v\in\mathcal{C}^{\infty}(M)$ and $t\in\mathbb{R}$ such that $\left|t\right|< \left\|v\right\|_{\infty}^{-1}$. Since $u$ is a minimum for (\ref{defalpha}),
\begin{multline} \label{inequ} \int_M \left|\nabla (u+tuv) \right|_g^2dv_g + \int_M (h+\alpha(g,h)) (u+tuv)^2dv_g \\- K_n^{-2} \left( \int_M (u+tuv)^{2^{\star}}dv_g \right)^{\frac{2}{2^{\star}}} \geq 0 \hskip.1cm. \end{multline}
Since $u$ satisfies (\ref{eullageq}), the left term in (\ref{inequ}) vanishes until the order $2$ in the Taylor development as $t\to 0$. Computing the second-order coefficient as $t\to 0$, one gets
\begin{multline} \label{ineq} \int_M \left|\nabla (uv) \right|_g^2dv_g + \int_M (h+\alpha(g,h)) (uv)^2dv_g - K_n^{-2}(2^{\star}-1) \int_M v^2 u^{2^{\star}}dv_g \\+ K_n^{-2}(2^{\star}-2)\left( \int_M v u^{2^{\star}}dv_g\right)^2 \geq 0 \hskip.1cm. \end{multline}
We now use the conformal transformation of the conformal Laplacian
\begin{equation} \label{translapla} \forall v\in\mathcal{C}^{\infty}(M), u^{2^{*}-1}\Delta_{\tilde{g}} v = \Delta_g(uv) - v\Delta_g u \end{equation}
where $\tilde{g} = u^{\frac{4}{n-2}}g$. We integrate (\ref{translapla}) against $uv$ and with (\ref{eullageq}),
\begin{eqnarray*}
\int_M \left|\nabla (uv)\right|_g^2dv_g &=& \int_M \left|\nabla v\right|_{\tilde{g}}^2 dv_{\tilde{g}} + \int_M v^2 u \Delta_g u dv_g \\
&=& \int_M \left|\nabla v\right|_{\tilde{g}}dv_{\tilde{g}}^2 - \int_M (h+\alpha(g,h))v^2u^2 dv_g + K_n^{-2} \int_M v^2 u^{2^{\star}}dv_g 
\end{eqnarray*}
and with (\ref{defKn}), (\ref{ineq}) becomes
\begin{equation} \label{hessianJ} \int_M \left|\nabla v\right|^2_{\tilde{g}} dv_{\tilde{g}} - n\omega_n^{\frac{2}{n}} \int_M \left(v - \int_M v dv_{\tilde{g}}\right)^2 dv_{\tilde{g}} \geq 0 \hskip.1cm. \end{equation}
This gives that $\lambda_1(\tilde{g}) \geq n\omega_n^{\frac{2}{n}}$. Note that if the inequality is strict for one solution $(h,u)$ of $F(u,\alpha(g,h),h)=0$, the theorem is proved.

We now assume that for any solution $(h,u)$ of $F(u,\alpha(g,h),h)=0$, we have $\lambda_1(u^{\frac{4}{n-2}}g) = n\omega_n^{\frac{2}{n}}$. We will apply the following theorem (\cite{HEN2005},Theorem 5.4,page 63) of Fredholm theory to $F$, with $U = W^{s+2,2}_+(M)\times \mathbb{R}$.
\begin{theo} \label{gene} Let $X$,$Y$ be two separable Banach spaces, $U$ an open set of $X$, $V$ a separable $\mathcal{C}^{\infty}$ Banach manifold and $F\in\mathcal{C}^{\infty}(U\times V,Y)$ which satisfy :
\begin{itemize}
\item For all $(u,v)\in F^{-1}(0)$, $DF(u)$ is surjective.
\item For all $(u,v)\in F^{-1}(0)$, $D_u F(u,v)$ is a Fredholm operator.
\end{itemize}
Then there exists a countable intersection of open dense sets (a residual set) $\Sigma \subset V$ such that for all $v\in\Sigma$, and for all $u\in F(.,v)^{-1}(0)$, $D_u F(u,v)$ is surjective. 
\end{theo}
Using (\ref{translapla}) and (\ref{defKn}), one gets for $(u,\beta,h) \in F^{-1}(0)$,
\begin{equation} \label{DiffFubeta} D_{(u,\beta)}F(u,\beta,h).(\theta,\mu) = u^{2^{*}-1} \left( \Delta_{\tilde{g}} \left(\frac{\theta}{u}\right) - n\omega_n^{\frac{2}{n}} \frac{\theta}{u} \right) + \mu u \end{equation}
where $\tilde{g}=u^{\frac{4}{n-2}}g$. Then, $D_{(u,\beta)} F(u,\beta,h)$ is a Fredholm operator. It remains to prove that if $(u,\beta,h) \in F^{-1}(0)$, $DF(u,\beta,h)$ is surjective. We have
\begin{equation} \label{DiffF} DF(u,\beta,h).(\theta,\mu,\tau) = u^{2^{*}-1} \left( \Delta_{\tilde{g}} \left(\frac{\theta}{u}\right) - n\omega_n^{\frac{2}{n}} \frac{\theta}{u} \right) + \mu u + \tau u \hskip.1cm. \end{equation}
$Im (D_{(u,\beta)}F(u,\beta,h))$ is a closed space in $W^{s,2}$ of finite codimension. Thus, since $Im (DF(u,\beta,h))$ contains $Im (D_{(u,\beta)}F(u,\beta,h))$, it is a closed space in $W^{s,2}$ by the following

\begin{lem} Let $X$ a banach space, and $E\subset F \subset X$ some subspaces. If $E$ is a closed finite co-dimentional subsbace of $X$, then $F$ is a closed subspace of $X$.
\end{lem}
\medskip {\bf Proof.}
Let $G$ a finite dimensional subspace of $X$ such that $X = E \oplus G$. We set $H = G\cap F$. Then, $F = E \oplus H$. Let $x_k \in F$ such that $x_k \to x$ as $k\to+\infty$. We denote $x_k = y_k + z_k$ with $y_k \in E$ and $z_k \in H$. 

We suppose that $(z_k)_{k\geq 0}$ is not bounded. Then, up to the extraction of a subsequence, $\left|z_k\right| \to +\infty$ as $k\to+\infty$. By Bolzano's theorem, up to the extraction of a subsequence, there exists $z\in H$ such that
$$ \frac{z_k}{\left|z_k\right|} \to z \hbox{ as } k\to+\infty \hskip.1cm.$$
Since $(x_k)$ converges as $k\to+\infty$,
$$ \frac{y_k}{\left|z_k\right|} = \frac{x_k}{\left|z_k\right|} - \frac{z_k}{\left|z_k\right|} \to -z \hbox{ as } k\to+\infty \hskip.1cm.$$
Since $E$ is closed, we get $z \in E\cap H = 0$, which contradicts $\left|z\right|=1$.

Then $(z_k)_{k\geq 0}$ is bounded and by Bolzano's theorem, up to the extraction of a subsequence, we can suppose that $z_k \to z \in H$ as $k\to +\infty$. Then,
$$ y_k = x_k - z_k \to x - z \hbox{ as } k\to+\infty \hskip.1cm.$$
and $y = x - z \in E$ since $E$ is closed. Therefore $x = y + z \in E+H = F$ and the proof of the lemma is complete.
\hfill $\diamondsuit$

\medskip Now, it suffices to prove that $Im (DF(u,\beta,h))^{\perp} = 0$, where $\perp$ refers to the orthogonal in $W^{s,2}$. Let $\phi \in Im (DF(u,\beta,h))^{\perp}$. Then, with (\ref{DiffF}),
$$ \forall \tau \in \mathcal{C}^m, \left\langle \phi , u \tau \right\rangle_{W^{s,2}} = 0 \hskip.1cm.$$
Since $u\in\mathcal{C}^m$ is positive and $\mathcal{C}^m$ is dense in $W^{s,2}$, we get $\phi = 0$.

\medskip By Theorem \ref{gene}, there exists $h \in V$ such that for all couple $(u,\beta)$ with $F(u,\beta,h)=0$, $DF_{(u,\beta)}(u,\beta,h)$ is surjective. We take in particular $\beta=\alpha(g,h)$ and we will deduce that for a minimal function $u$, $\lambda_1(\tilde{g}) = n\omega_n^{\frac{2}{n}}$ is simple with $\tilde{g} = u^{\frac{4}{n-2}}g$. We claim that
\begin{equation}\label{propphi} \forall \phi\in E_1(\tilde{g})\setminus\{0\}, \int_M u^2 \phi dv_g \neq 0 \hskip.1cm. \end{equation}
Indeed, if $\phi$ is an eigenfunction for $\lambda_1(\tilde{g})$ such that this integral vanishes, one easily checks with (\ref{DiffFubeta}) that $u \phi$ is orthogonal to the image of $D_{(u,\beta)}F(u,\alpha(h,g),h)$ in $L^2(g)$. It implies $\phi=0$ and we obtain (\ref{propphi}). Since a bounded linear form vanishes on a one-codimensional space, we get that $\lambda_1(\tilde{g})$ is simple. Thus, $\lambda_1(\tilde{g})$ cannot be an extremal eigenvalue in the sense of \cite{ELS2007} and as a result, $\lambda_1(\tilde{g}) = n\omega_n^{\frac{2}{n}}$ is not locally maximal. The proof of Theorem \ref{rig} for $n\geq 3$ is complete.

\medskip\textbf{Acknowledgements}

I would like to thank my thesis advisor O. Druet for having pointed out to me the interest of Sobolev inequalities for the existence of metrics with large first eigenvalue and the referee for remarks on the manuscrit which led to improvements in the presentation.

\bibliographystyle{plain}
\bibliography{biblio-rigidity}

\begin{thebibliography}{10}

\bibitem{AUB1976}
T.~Aubin.
\newblock \'{E}quations diff\'erentielles non lin\'eaires et probl\`eme de
  {Y}amabe concernant la courbure scalaire.
\newblock {\em J. Math. Pures Appl. (9)}, 55(3):269--296, 1976.

\bibitem{COL2004}
B.~Colbois.
\newblock Spectre conforme et m\'etriques extr\'emales.
\newblock In {\em S\'eminaire de {T}h\'eorie {S}pectrale et {G}\'eom\'etrie.
  {V}ol. 22. {A}nn\'ee 2003--2004}, volume~22 of {\em S\'emin. Th\'eor. Spectr.
  G\'eom.}, pages 93--101. Univ. Grenoble I, Saint, 2004.

\bibitem{COL1994}
B.~Colbois and J.~Dodziuk.
\newblock Riemannian metrics with large {$\lambda_1$}.
\newblock {\em Proc. Amer. Math. Soc.}, 122(3):905--906, 1994.

\bibitem{COL2003}
B.~Colbois and A.~El~Soufi.
\newblock Extremal eigenvalues of the {L}aplacian in a conformal class of
  metrics: the `conformal spectrum'.
\newblock {\em Ann. Global Anal. Geom.}, 24(4):337--349, 2003.

\bibitem{DRU2002}
O.~Druet.
\newblock Optimal {S}obolev inequalities and extremal functions. {T}he
  three-dimensional case.
\newblock {\em Indiana Univ. Math. J.}, 51(1):69--88, 2002.

\bibitem{DRH2002}
O.~Druet and E.~Hebey.
\newblock The {$AB$} program in geometric analysis: sharp {S}obolev
  inequalities and related problems.
\newblock {\em Mem. Amer. Math. Soc.}, 160(761):viii+98, 2002.

\bibitem{ELS1986}
A.~El~Soufi and S.~Ilias.
\newblock Immersions minimales, premi\`ere valeur propre du laplacien et volume
  conforme.
\newblock {\em Math. Ann.}, 275(2):257--267, 1986.

\bibitem{ELS2007}
A.~El~Soufi and S.~Ilias.
\newblock Laplacian eigenvalue functionals and metric deformations on compact
  manifolds.
\newblock {\em J. Geom. Phys.}, 58(1):89--104, 2008.

\bibitem{HEB1992}
E.~Hebey and M.~Vaugon.
\newblock Meilleures constantes dans le th\'eor\`eme d'inclusion de {S}obolev
  et multiplicit\'e pour les probl\`emes de {N}irenberg et {Y}amabe.
\newblock {\em Indiana Univ. Math. J.}, 41(2):377--407, 1992.

\bibitem{HEN2005}
D.~Henry.
\newblock {\em Perturbation of the boundary in boundary-value problems of
  partial differential equations}, volume 318 of {\em London Mathematical
  Society Lecture Note Series}.
\newblock Cambridge University Press, Cambridge, 2005.
\newblock With editorial assistance from Jack Hale and Ant{\^o}nio Luiz
  Pereira.

\bibitem{KON2010}
G.~Kokarev and N.~Nadirashvili.
\newblock On first {N}eumann eigenvalue bounds for conformal metrics.
\newblock In {\em Around the research of {V}ladimir {M}az'ya. {II}}, volume~12
  of {\em Int. Math. Ser. (N. Y.)}, pages 229--238. Springer, New York, 2010.

\bibitem{LED2000}
M.~Ledoux.
\newblock The geometry of {M}arkov diffusion generators.
\newblock {\em Ann. Fac. Sci. Toulouse Math. (6)}, 9(2):305--366, 2000.
\newblock Probability theory.

\bibitem{LIY1982}
P.~Li and S.~T. Yau.
\newblock A new conformal invariant and its applications to the {W}illmore
  conjecture and the first eigenvalue of compact surfaces.
\newblock {\em Invent. Math.}, 69(2):269--291, 1982.

\bibitem{PETR2013}
R.~Petrides.
\newblock Existence and regularity of maximal metrics for the first {L}aplace
  eigenvalue on surfaces.
\newblock 2013.

\bibitem{SCH1984}
R.~Schoen.
\newblock Conformal deformation of a {R}iemannian metric to constant scalar
  curvature.
\newblock {\em J. Differential Geom.}, 20(2):479--495, 1984.

\bibitem{YAN1980}
P.~C. Yang and S.~T. Yau.
\newblock Eigenvalues of the {L}aplacian of compact {R}iemann surfaces and
  minimal submanifolds.
\newblock {\em Ann. Scuola Norm. Sup. Pisa Cl. Sci. (4)}, 7(1):55--63, 1980.

\end{thebibliography}

\end{document}